\documentstyle [12pt,amsfonts,amssymb]{article}

\newtheorem{lemm}{Lemma}      
\newtheorem{prop}{Proposition}   
\newtheorem{theo}{Theorem}   

\newtheorem{defi}{Definition}

\newenvironment{proo}{\noindent {\bf Proof:}}{$\Box$ \vspace{7mm}}

\DeclareMathAlphabet{\scri}{OT1}{eus}{m}{n}

\newcommand{\labl}[1]{\label {#1}}

\newcommand {\text}{\mbox}

\newcommand {\dom} {\text{\rm{dom}}}

\newcommand {\fld} {\text{\rm{fld}}}
\newcommand {\defin} {\text{\rm{def}}}

\begin{document}

\title{The Consistency of ZFC$ \,+\, 2^{\aleph _{0}} > \aleph _{ \omega } + {\cal I } (\aleph _{2} )
 = {\cal I }(\aleph _{\omega } ) $}
\author{Martin Gilchrist and Saharon Shelah}
\maketitle

\section{Introduction}

The  \footnote{S. Shelah partially supported by
a research grant from the basic research fund of the Israel Academy of Science;
Pul. Nu. 583.}
basic notion that will be studied in this work is than of an {\em identity}.   It arises naturally in
a Ramsey theory setting when considering the coloring patters on finite sets
that occur when coloring infinite complete graphs with infinitely
many colors.  We first give some definitions and establish some notation. 

An {\em $\omega$-coloring} is a pair $ \langle f,B \rangle$ where $f:[B]^{2} \longrightarrow \omega$.  The set $B$ is the {\em field} of $f$ and denoted Fld$(f)$.
\begin{defi} Let $f,g$ be 
$\omega$-colorings.  We say that $f${\em realizes} the coloring $g$ if
there is a one-one function $ k : \fld (g) \longrightarrow
\fld  (f) $ such that for all $ \{ x,y\} ,\{ u,v \} \in \dom (g)$
$$ f(\{ k(x), k(y) \} ) \not=  f(\{ k(u), k(v) \} ) \Rightarrow g(\{ x,y \} )
 \not=  g(\{ u,v \} ).$$  We write $f \simeq g$ if $f$ realizes $g$
and $g$ realizes $f$.  It should be clear that $ \simeq$ induces an
equivalence relation on the class of $\omega$-colorings.  We call the 
$ \simeq$-classes of $ \omega$-colorings with finite fields {\em identities}.
\end{defi}

If $ f,g,h,k$ are $ \omega $-colorings, with $ f \simeq g $ and $ h \simeq k $,
then $f$ realizes $h$ if and only if $g$ realizes $k$.  Thus without risk
of confusion we may speak of identites realizing colorings and of identities
realizing other identities.  We say that an identity $I$ is of {\em size r }
if $| \fld (f) | = r $ for some (all) $f \in I $.

Let $\kappa $ be a cardinal and $ f : [\kappa]^{2} \longrightarrow \omega$.
We define ${\cal I}(f)$ to be the collection of identities
realized by $f$ and ${\cal I }  (\kappa)$ to be $\bigcap \{ {\cal I }  (f)|
f: [\kappa] ^{2} \longrightarrow \omega\}$.  We now define a specific 
collection of identities.  Let $h : ^{< \omega }\!\!2 \longrightarrow \omega$
be 1-1.  Define $ f : [2^{\omega } ] ^{2} \longrightarrow \omega$ by
$f(\{ \alpha, \beta \} )= h ( \alpha \bigcap \beta )$. We define
${\cal J } = {\cal I } (f)$.  Note that ${\cal J }$ is independent of the choice of $h$.  In \cite{Shtwocard}, the second author proved
that $2 ^{\aleph _{0}}
> \aleph _{\omega }$ implies ${\cal I } (\aleph _{\omega}) = {\cal J }$.

In \cite{GSh1}, was shown consistentcy of $ZFC + {\cal I } (\aleph _{2} ) \not= {\cal I }(\aleph _{\omega})$.  Here we will show 

\medskip 

\noindent {\bf Main Theorem.}   If $ZFC$ is consistent then $ZFC + 2^{\aleph _{0}} > \aleph _{ \omega } + {\cal I } (\aleph _{2} )
 = {\cal I }(\aleph _{\omega } )$ is
consistent.
  
\medskip

This is accomplished by adding $ \nu > \aleph _{\omega}$
random reals to a model of $GCH$.  As  $2 ^{\aleph _{0}} > \aleph _{\omega }$ holds in the resulting model we need only show that ${\cal I }(\aleph _{2})
\supseteq{\cal J }$ is true.

\medskip

\section{The Partial Order}
We establish the notation necessary to add many random reals to a model
of ZFC.  For a more detailed explanation see \cite{Jech}.  Let $ \nu > \aleph _{\omega}$ be a cardinal.  Let $ \Omega = \,\,^{\nu } \! \{ 0,1 \}.$  Let $T$
be the set of functions $t$ from a finite subset of $ \nu $ into $ \{ 0,1 \} $.
For each $ t \in T$, let $ S _{t} = \{ f \in \Omega : t \subset f \}$ and let
$ {\cal S } $ be the $ \sigma$-algebra generated by $ \{ S _{t} : t \in T \}$.
The product measure $m$ on ${\cal S }$ is the unique measure so that
$ m ( S _{t} ) = \frac{ 1} { 2 ^{|t|}}$.  We define $ {\cal B } _{1}$ to be
the boolean algebra $ {\cal S } / J $ where $J$ is the ideal of all
$ X \in {\cal S} $ of measure $0$.  We define a partial order
$\langle {\Bbb P } , < \rangle $ by letting
$ {\Bbb P } = {\cal B } _{1}\setminus J $ and the order be
inclusion modulo $J$.  The following two theorems can be found in \cite{Jech}.

\begin{theo}  ${\Bbb P }$ is c.c.c. 
\end{theo}

\begin{theo} Let $M$ be a model of set theory and $G$ be ${\Bbb P }$-generic.
Then $M[G]$ satisfies $ 2^{\aleph _{0} } \geq \nu$.
\end{theo}
  Let $ Y = \{ y _{\alpha } : \alpha < \nu \}$.  Let $ \Gamma$ denote the collection of all $ \tau (\bar{y})$
where $ \bar{y}$ is a tuple from $Y$ and $ \tau (\bar{x} )$ is a boolean
term with free variables $ \bar{x}$.  For $ \alpha < \nu$ denote by
$ t_{\alpha} \in T $ the function whose domain is
$ \{ \alpha \}$ such that $ t_{\alpha} ( \alpha) = 0$.  There is an obvious embedding of $ \Gamma $ into $ {\cal S }$
which extends the map $ y _{\alpha } \mapsto S_{t _{\alpha} }$ and respects the boolean operations.
We denote by ${ \cal B } _{0} $ the image of $ \Gamma$ in ${\cal S }$.  It
should be clear that $ {\cal B } _{0} $ is a boolean algebra.  We call the
elements of $Y$ {\em generators}.  Elements of
${\cal B }_{0}$ are denoted by their preimage in $\Gamma$.  The following 
theorem should be clear.

\begin{theo}\labl{theo33}  For $p\in {\cal S } $ and $ \epsilon > 0 $ there exists
a finite $ u \subset Y $ and a boolean formula $ \tau (\bar {x})$ such
that $ \mu ( \tau ( \bar {u} ) \triangle p) < \epsilon$, where $ \triangle$
denotes the symmetric difference.
\end{theo}

\section{A Combinatorial Statement}
Here we formulate a combinatorial statement $ [I, \kappa,\lambda,g,f]$ which will
play a crucial role in the proof of the main result.  We require some preliminary definitions.  Let $Y,{\cal S }, {\cal B}_{0}, {\cal B} _{1}, \mu $ and $
{\Bbb P } $  be as in the previous section.  Let $g, f : \omega \longrightarrow \omega$.
 For each
$L < \omega $ let $ {\cal T } _{L} $ be a finite set of boolean
terms $\tau(\bar{x})$ where $ \bar{x}= ( x _{1} , \ldots ,x_{f(L) } )$
which is complete in the sense that for any boolean term $ \sigma (\bar{x})$
there is some $ \tau(\bar{x}) \in {\cal T}_{L}$ such that $\sigma(\bar{x})=
\tau(\bar{x})$ is a valid formula of the theory of boolean algebras.  Let $ {\cal T } = \bigcup \{ {\cal T } _{L} : L< \omega \}$. In the following we work only with boolean formulas in $ {\cal T } $. List $ {\cal T }_{L}$ as 
$ \{ \tau ^{L} _{i} : i \leq h(L) \}$. 
\medskip
For $ L < \omega $ define $ {\Bbb T}_{L}= ( {\cal T} _{L}) ^{g(L)}.$  
For $ w \in [\kappa]^{2} $ and $ L < \omega $ define $${\Bbb T} _{w,L}=
\{ \langle \tau _{1} ( \bar{x} ^{w,t}_{L}), \ldots , \tau _{g(L)} ( \bar{x}
^{w,t }_{L} ) \rangle : t= \langle \tau _{1}, \ldots, \tau _{g(L)} \rangle
\in {\Bbb T } _{L} \}$$  where $\bar{x}^{w,t}_{L}=
\langle x^{w,t} _{L,1}, \ldots , x ^{w,t} _{L,f(L)} \rangle$ is a sequence of distinct variables for each triple $ ( w,t,L ), $ and where $$ \bar{x} ^{ w,t }_{L} \cap
\bar{x} ^{ v,u}_{M} \not= \emptyset \Rightarrow ( t=u \wedge w=v \wedge
L=M ).$$  Let  $X$ denote $$ \bigcup \{ \bar{x}^{w,t}_{L} : t \in {\Bbb T}_{L},\, L < \omega,\,
w \in [\kappa ] ^{2} \}.$$  Let $ {\cal C} (P,L)$ denote $$\{ c : c {\text{ is a
mapping of } } [P]^{2} {\text{ into }} \{ 1, \ldots, g(L) \} \}.$$

\begin{defi}Let $ k,m < \omega $ and $\langle \tau_{n} ( \bar{x}) : n \leq k \rangle$ be a sequence of $m$-ary boolean formulas.  Let $\bar{u}$ be an $m$-tuple from
$ Y$.  Then $\langle \tau _{n}( \bar{u}) : n \leq k \rangle $ is called a {\em 
partition sequence} if $ \mu (\tau_{m} ( \bar{u}) \cap  \tau_{n} ( \bar{u})) = 0$ for all $ m,n $ with $m \not= n$, and $ \mu ( \bigcup \{ \tau _{n} ( \bar{u} ) 
: n \leq k \} ) = 1 $.
\end{defi}

The combinatorial statement will now be defined.  
\begin{defi}Let $ I $ be an $r$-identity, $ \lambda \leq \omega$ and $ \kappa $ a
cardinal.  We say that $ [I,\kappa, \lambda,g,f] $ {\em holds} if the following is true:
there exist $ \bar{u} _{w,L},\, \tau ^{w} _{L,m} ( w \in [\kappa] ^{2} ,\, L < 
\lambda,\, 1 \leq m \leq g(L))$ such that for all $ w \in [\kappa] ^{2}, \,L< 
\lambda$ and $ P \in [\kappa] ^{r}$
\begin{enumerate}
\item [C1.]$\bar{u} _{w,L}$ is a tuple in $Y$ of length $ f(L)$
\item [C2.] $ \tau ^{w} _{L,m } \in {\cal T } _{L}, \langle \tau ^{w} _{L,1}
, \ldots , \tau ^{w} _{L,g(L)} \rangle \in {\Bbb T} _{L}$
\item [C3.]$ \langle \tau ^{w} _{L,m } (\bar{u} _{w,L} ) : 1 \leq m \leq g(L) \rangle$ is
a partition sequence
\item [C4.]for $ N \leq L  , \mu ( \bigcup  \{ \tau^{w} _{ N,m} (\bar { u} _{w,N}) \cap \tau^{w} _{  L,m} (\bar { u} _{w,L})  : m \leq g(N) \} )\geq 1 - 1/2^{N} $
\item [C5.] the measure of 
$$ \bigcup \{  \bigcap  \{\tau ^{z}_{L,c(z)}(\bar{u} _{z,L}) : z \in [P]^{2} \} : c \in 
{\cal C}(P,L) \wedge c {\text{ realizes } } I \} $$ is less that $ 1/L$.
\end{enumerate}
\end{defi}

\section{Proof of the Main Theorem}
The theorem follows from the following three lemmas which will be proved later.
\begin{lemm}\labl{lemm1} Let $I \in {\cal J}$.  For no $g, f:\omega \longrightarrow \omega $ and $\kappa > \aleph _{\omega}$ do we have $ [I, \kappa, \omega,g,f]$.
\end{lemm}

\begin{lemm}\labl{lemm2} Let $I \in {\cal J } $, $ \kappa  \geq \aleph _{0}$ and $
g,f:\omega \longrightarrow \omega$ be such that $ [I, \kappa, \omega ,g,f]$ fails.
Then there exists $m< \omega $ such that $[I, m,m,g,f]$ fails.
\end{lemm}

\begin{lemm}\labl{lemm3} Let $I \in {\cal J }$ and $M$ be a model of set theory satisfying GCH.  Let $G$ be ${\Bbb P}$-generic over $M$.  If it is true in
$M[G]$ that $ I \not\in {\cal I}(\aleph _{2})$, then in $M$ there exists
$g,f:\omega \longrightarrow \omega$ such that $[I,m,m,g,f]$ holds for all
$m < \omega$. 

\end{lemm}
 We suppose that these lemmas are true and prove the main result.  Let $M$ be a model of $ ZFC + GCH$.  Let $I \in {\cal J } $ and towards a contradiction 
suppose that $ I \not\in {\cal I} ( \aleph _{2} )$ in $M[G]$ where $G$
in $ {\Bbb P } $-generic over $M$.  By lemma \ref{lemm3} in $M$ there exist
$g,f : \omega \longrightarrow \omega $ such that $[I,m,m,g,f]$ holds for all
$m < \omega$.  But from lemma \ref{lemm1},
$ [ I,(\aleph _{\omega} )^{+},\omega,g,f]$ fails, and so by lemma \ref{lemm2}
there exists $ m < \omega$ such that $ [I,m,m,g,f]$ fails, contradiction.

\subsection{Proof of the first lemma }
Assume that the conclusion of lemma  fails.  Let $ \kappa > \aleph_ {\omega}$.  Let $g, f : \omega \longrightarrow \omega$
be such that $[I,\kappa, \omega,g,f] $ holds.  We force
with the partial order $ {\Bbb P}$, where ${\Bbb P }$ is defined with $ \nu = \kappa $.  Let $G \subseteq {\Bbb P } $ be a generic set. 
For $ L < \omega$ we define $ c_{L} : [\kappa] ^{2} \longrightarrow \omega$
by $ c_{L} (w ) = m $ if $ \tau^{w}_{L,m} (\bar{u}
_{w, L } )/J \in G$.

\begin{prop}  For all $ w \in [\kappa] ^{2}$ there exists $ N< \omega, m< \omega $ such that $c_{L} ( w) = m$ for all 
$ L > N$.
\end{prop}

\begin{proo}  For $ w \in [\kappa] ^{2} $ define $$D_{w}=
\{ p \in {\Bbb P } : p \Vdash \exists N \exists m ( c_{L} ( w) = m 
{\text{ for all }} L > N \}.$$  We claim that $ D_{w } $ is dense 
in $ {\Bbb P}$.  To this end choose $ p^{*} \in {\Bbb P}$ and let $ p \in {\cal
S} $ be such that $ p / J = p^{*}$.  Let $ \mu ( p) =
\delta$.  As $ \delta > 0$ we
can choose $N$ such that $\sum _{ L > N} 1/2^{L} < \delta / 3.$  By
C4 of
the definition of $ [I, \kappa, \omega ,g,f],$ 
$$\mu (\bigcup\{ \bigcap \{ \tau^{w} _{ L ,m } ( \bar{u} _{w, L} ): L > N \} : m \leq g(N) \}  )
> 1 - (\delta / 3).$$
Thus  $$\mu (\bigcup\{ \bigcap \{ \tau^{w} _{ L ,m } ( \bar{u} _{w, L} ): L> N \} : m \leq g(N) \} \cap p )
> \delta / 3.$$
There is thus an $m  \leq g(N) $ such that $ \mu (q) > 0$, where
 $$ q= \bigcap  \{ \tau ^{w}_{L ,m } ( \bar{u} _{w, L} ):
L > N \}\cap p.$$  Clearly $q / J \Vdash c_{L} ( w ) =m $
for all $ L > N$.  Thus the proposition is proved. 
\end{proo}

We now continue with the proof of the lemma. Define $ c : [\kappa]^{2} \longrightarrow \omega$ in $ M[G] $ by $ c(w
) = \lim _{L \longrightarrow \omega } c_{L} ( w )$.  Fix $ P \in [\kappa]^{r}$.
By property C5 of $[I, \kappa, \omega,g,f],$
$$\sup \{\mu (  p ): p / J \Vdash 
``c_{L} {\text{ realizes }}I {\text{ on }} P" \} <1/ L.$$
Thus  $$\sup  \{\mu( p) : p/J \Vdash 
``c {\text{ realizes }} I {\text{ on }} P" \}  < 1/ L$$ for all sufficiently large
$ L < \omega $.  Hence this set has measure 0 and so it is true
that $c$ does not induce $I$ in any generic extension.  A contradiction occurs as
$ \kappa > \aleph _{ \omega}$ and by \cite{Shtwocard} every coloring $ c: [
\kappa] ^{2} 
\longrightarrow \omega $ must realize $I$.  Thus the lemma is proved.

\subsection{ Proof of the second lemma  }
The proof of lemma \ref{lemm2} is accomplished by showing that it
is possible to represent the statement $ [I,\kappa, \omega, g,f]$ by a theory in a language of propositional constants when the propositional constants are 
assigned suitable meanings.  The compactness theorm is then
used to show that the failure of $[I ,\kappa, \omega,g,f] $ implies the failure
of $ [I,m,m,g,f]$ for all sufficiently large $m$ in $ \omega$.

\medskip
Throughout this section fix $g, f : \omega \longrightarrow \omega $.  Let $ 
{\cal B }_{0} $ and
$ \mu $ be as previously defined.  Let $I$ be an $r$-identity for some $ r < \omega$. Consider $X$, the collection of free variables previously
defined.  Define ${\cal L}  = \{
p_{w} : w \in [X] ^{2} \}$ to be a collection of propositional
constants.  For each partition $ {\cal P }$ of $X$ let  $\sim _{{\cal P }}$ denote
the associated equivalence relation. 
Let $${\cal A}: [\kappa]^{2} \times \{ (L,m) : L < \omega \wedge 1 \leq m
\leq g(L) \} \longrightarrow {\cal T}$$ be such that ${\cal A } ( w, L,m) \in
{\cal T} _{L}$ for all $ w \in [\kappa] ^{2} $ and $ 1 \leq m \leq g(L)$.
Let $${\cal Q} = \{ q^ {w} _{L,m,i } : w \in [\kappa] ^{2} , L < \omega , 1 \leq m \leq g(L), i \leq h(L) \}$$ be a collection of propostional constants. 
Denote $ {\cal R} = {\cal L} \bigcup {\cal Q}$.
For each ${\cal P } $ a partition of $X$ and function $ {\cal A}$ define a
truth valuation $ V_{{\cal P }, {\cal A} }: {\cal R} \longrightarrow
\{\bf{T}, \bf{F} \} $ by $ V_{ {\cal P}, {\cal A}} ( p_{w} ) = \bf{T} $ iff $ w = \{ i,j\} \wedge
i \sim _{{\cal P }} j$ and $V_{ {\cal P}, {\cal A}} (  q^ {w} _{L,m,i }) = \bf{T}$ iff
$ {\cal A} ( w,L,m) = \tau ^{L} _{i}$.  There is a propositional theory
 $T_{0}$ such that a truth valuation $V$ models $T_{0}$ if and only if
$V = V _{{\cal P }, {\cal A } } $ for some function ${\cal A} $ and partition
${\cal P }$.

Let $ V  $ be a truth valuation that models the
theory $T_{0}$.  Denote by ${\cal P}_{V}$ the partition of $X$ defined by $ x _{1}  \sim_{{\cal P }_{V}}
x_{2}  \Leftrightarrow V ( p _{\{ x_{1}, x_{2} \}} ) = \bf{T}$.  Fix a mapping
$ v_{V} : X \longrightarrow Y$ such that $ v_{V} (x) = v_{V} (y)
\Leftrightarrow x \sim _{ {\cal P } _{V}} y $.  For $ L < \omega , 1 \leq m \leq g(L)
$ and $ w \in [ \kappa ]^{2} $ define $ \tau ^{V, w} _{ L, m } $ to be
$ \tau ^{L} _{i}$ if $ V ( q ^{w} _{L,m,i} ) = \bf{T}$.  Let $t =  t ^{V,w} _{L}
$ denote $\langle \tau ^{V,w} _{L,1}, \ldots , \tau ^{V,w} _{L, g(L)} \rangle \in {\Bbb
T} _{L}$. For each such sequence let $ \bar {x} ^{V,w,t} _{L} $ denote
$ \bar{x} ^{w,t} _{L}$ and write
$ \tau ^{V,w} _{L,m} (\bar{u}^{V,w} _{L})$ for the ${\cal B} _{0}$-term
obtained from $ \tau ^{V,w} _{L,m}(\bar{x} ^{V,w,t} _{L}) $ by substituting the variables
$ \bar{x} ^{V,w,t} _{L} $ by their image under $ v_{V}$. Note that since $ {\Bbb T}_{L}$
is finite, for each $ L < \omega $ and $ w \in [\kappa]^{2}$,
$$ X_{L}^{w} = _{\defin} \bigcup \{\bar{x} ^{V,w,t} _{L}
:t =  t ^{V,w} _{L} \in {\Bbb T} _{L} \wedge
 V {\text{ models }} T _{0} \}$$ is finite.

\begin{lemm}

Let $ k < \omega $ and $ \sigma ( x_{1}, \ldots, x _{k}) $ be a boolean term.
For $ 1 \leq i \leq k $ let $ L _{i} < \omega, 1 \leq m_{i} \leq g(L_{i})$ and $
w _{i} \in [ \kappa ] ^{2}$.  Let $ \theta(y)$ be a statement of one of the
 forms $ \mu (y) < 1/n , \mu(y) > 1/n $ or $ \mu(y)
=0$, where $y$ runs through $ {\cal B } _{0}$.  There exists a propositional formula $ \chi $ such that for all valuations
$V$ modelling $T_{0}, V $ models $ \chi $ if and only if 
$ \theta (\sigma (\tau ^{V, w_{1}} _{L_{1}, m_{1} } (\bar{u} ^{V,w_{1}} _{L_{1}})
, \ldots, \tau ^{V, w_{k}} _{L_{k}, m_{k} } (\bar{u} ^{V,w_{k}} _{L_{k}}))).$
\end{lemm}
\begin{proo} Let $ W = \bigcup \{ X^{w_{i}} _{L_{i}}: 1 \leq i \leq k \}$. Define $ {\cal V } = \{V:
V $ is a truth valuation modelling $ T_{0} \}$.  Since $ {\cal T } _{L_{i}}$ is finite for all $ 1 \leq i \leq k $ the
collection $ S = \{  \langle \tau ^{V, w _{i}} _{L_{i}, m_{i}} : 1 \leq i \leq k
\rangle : V \in {\cal V } \}$ is a finite set.  For each $ s  \in S $ define 
$ {\cal V } _{s} = \{ V \in {\cal V } :\langle  \tau
^{V,w_{i}}_{L_{i}}  : 1 \leq i \leq k
\rangle = s \}$

For the moment fix $ s \in S $.  Each $ V \in {\cal V }_{s} $
induces a partition, $ {\cal P } _{V_{s}} $ of $X$ and thus of $W$. Since every permutation of $Y$ induces an automorphism of $ {\cal B } _{0}$ which preserves
the measure,
for $ V _{1}, V _{2} \in {\cal V _{s}}, {\cal P} _{V_{1}} \restriction W = 
{\cal P } _{V _{2}} \restriction W$ implies $$\mu ( \sigma (\tau ^{V_{1}, w_{1}} _{L_{1}, m_{1} } (\bar{u} ^{V_{1},w_{1}} _{L_{1}})
, \ldots, \tau ^{V_{1}, w_{k}} _{L_{k}, m_{k} } (\bar{u} ^{V_{1},w_{k}} _{L_{k}})))$$ $$= \mu ( \sigma (\tau ^{V_{2}, w_{1}} _{L_{1}, m_{1} } (\bar{u} ^{V_{2},w_{1}} _{L_{1}})
, \ldots, \tau ^{V_{2}, w_{k}} _{L_{k}, m_{k} } (\bar{u} ^{V_{2},w_{k}} _{L_{k}}))).$$  As there are only finitely many partitions of $ W $ there
is a formula $ \chi_{s} $ that chooses those partitions in $\{ {\cal P} _{V}
: V \in {\cal V_{s} } \} $ that produce the desired measure.
We define $ \chi = \bigvee _{s \in S } (\eta _{s} \Rightarrow \chi _{s})$,
where $ \eta _{s} $ is a formula such that $ V \in {\cal V} $ implies
$s = \langle \tau ^{V,w_{i}} _{L_{i},m_{i} } : 1 \leq i \leq k \rangle$ if and only if $ V(\eta _{s}) = \bf{T}$.
\end{proo}

\begin{lemm}\labl{lemm4.6}  There is a propositional theory $T$ such that $T$ is consistent
if and only if $[I, \kappa, \omega, g,f ]$ holds.
\end{lemm}
\begin{proo}  By the previous lemma, for each triple $(w,L, P)$ where $ w \in [\kappa]^{2}, L < \omega$ and $P \in [\kappa]^{r}$ there exists a formula
$ \chi_{w,L,P}$ such that a truth valuation $V$ models $ T_{0}
\bigcup \{ \chi _{w,L,P} \}$ implies C1-C5 hold for $ w, L, P$ and the sequences of
boolean terms and generators defined by the valuation.  We define
$T$ to be $ T_{0} \bigcup \{ \chi_{w,L,P} :  w \in [\kappa]^{2}, L < \omega$ and $P \in [\kappa]^{r} \}$.  It is easily seen that the consistency of $T$ 
implies that
$ [ I, \kappa, \omega, g ,f ] $ holds.  In this regard one should observe that $Y$ is large enough to realize any desired partition.

Now suppose that $ [ I, \kappa, \omega, g ,f ] $ holds.  The existence of the sequences of terms $ t^{w} _{L} = \langle \tau ^{w} _{L,1}, \ldots,
\tau ^{w} _{L, g(L)} \rangle $ and generators $\bar{u} _{w,L} =
\langle u _{w,L,1} , \ldots, u _{w,L, f(L)} \rangle $ defines a function
$ {\cal A }$ and partition $ {\cal P}$ in the following manner.
Let ${\cal A} ( w ,L,m) = \tau ^{L} _{i} $ if $ \tau ^{w} _{L,m} = \tau ^{L} _{i}$.
A partition $ {\cal P } ^{\prime}$ of $ \bigcup\{ \bar{x} ^{w,t} _{L}
: t = t^{w} _{L}, w \in [\kappa]^{2}, L < \omega \}$ is first defined
by setting $ x ^{w,t} _{L,i} \sim _{{\cal P } ^{\prime}} x^{v,u} _{M,j}$ if
$ u_{w,L,i} = u_{v,M,j}$ where $ t= t^{w} _{L} $ and $ s = t^{v} _{M}$.  We choose
a partition of $X$ which is an extension of $ {\cal P } ^{\prime} $ and denote
it by $ {\cal P }$.  The truth valuation
$ V _{ {\cal P } , {\cal A} } $ models the theory $T$.  This completes
the proof of lemma \ref{lemm4.6}.
\end{proo} 

Lemma \ref{lemm2} follows from the compactness theorem for propositional
logic.

\subsection {Proof of the third lemma  }
Towards a contradiction let $I $ be an identity on $r < \omega $ elements, $d$ a $ {\Bbb P}$
-name for a function and $p \in {\Bbb P} $ such that 
$$ p \Vdash `` d:[
\aleph_{2}] ^{2} \longrightarrow \omega \wedge d {\text{ does not realize }}I ".$$  Without loss
of generality we assume that $ p = 1 _ {\Bbb P }$. 
For each $w\in[ \aleph_{2}]^{2} $ choose a sequence 
$ \langle b_{ n }^{w} : n < \omega \rangle $ and a sequence 
$\langle p^{w} _{ n } : n < \omega \rangle \in [{\cal S}]^{\omega}$ such that $\langle p^{w}_{n}/J : n < \omega \rangle$ is a maximal antichain
in $ {\Bbb P}$ and $ p^{w} _{n}/ J 
\Vdash d(w) = b^{w} _{n}$.   Let $
b: [\aleph _{2}]^{2} \times \omega \longrightarrow \omega$
be defined by $b(w,n) = b^{w} _{n}$. 

For $ w \in[ \aleph _{2}]^{2} , L < \omega $ choose $ g(w,
L ) $ so that $ \sum _{n > g(w,L ) } \mu ( p_{n}^{w}) < 
1/(2^{L+5} L)$. The next lemma follows from theorem \ref{theo33}.

\begin{lemm}\labl{lemm4.19} There exists a function $f : [\aleph_{2}]^{2} \times \omega \longrightarrow \omega$ sequences of boolean terms $ \langle \sigma _{L,m}^{w}: m \leq g(w,L) \rangle$ and generators
$\bar{v} _{w,L}  ( w \in [\aleph_{2}]^{2} , L < \omega)$ such that:

\begin{enumerate}
\item $ \bar{v} _{ w,L } = \{ y_{w,L,k} : k \leq h(w,L) \}$
\item For $ m \leq g(w,L)$ we have $$ \mu ( p_{m}^{w} \triangle \sigma _{L,m}^{w} (\bar{v} _{w,L})) < \frac{ 1}{(L2^{L+5} [g(w,L)] ^ {r^{2}+1} )}.$$
\end{enumerate}
\end{lemm}

\begin{lemm}\labl{lemm99} There exists a function $f : [\aleph_{2}]^{2} \times \omega \longrightarrow \omega$ sequences of boolean terms $ \langle \rho _{L,m}^{w}: m \leq g(w,L) \rangle$ and generators
$   \bar{v} _{w,L}( w \in [\aleph_{2}]^{2} , L < \omega)$ such that:
\end{lemm}
\begin{enumerate}
\item $ \bar{v} _{ w,L } = \{ y_{ w,L,k} : k \leq f(w,L) \}$
\item $  \langle \rho _{L,m}^{w} (\bar{v} _{w,L})
: m \leq g(w,L) \rangle $ is a partition sequence
\item  For $ m < g(w,L)$ we have $$ \mu (   p_{m}^{w} \triangle \rho _{L,m}^{w} (\bar{v} _{w,L})) < \frac{1} {2^{L+3} L[g(w,L)] ^ {r^{2}}  }$$
\item $ \mu (
 p_{g(w,L)}^{w} \triangle \rho _{L,g(w,L)}^{w} (\bar{v} _{w,L})) < 
\frac {1}{  L2^{L+3} }$.
 
\end{enumerate}

\begin{proo}  Let $ f, \sigma ^{w} _{L,m}, $ and $ \bar{v} _{w,L} $ satisfy
the conclusion of the last lemma.  For $ m < g(w,L) $ define
$\rho _{L,m} ^{w}(\bar{v} _{w,L}) =\sigma _{L,m}^{w} (\bar{v} _{w,L})\setminus \bigcup \{\sigma _{L,i}^{w} (\bar{v} _{w,L}) : i < m \} $.  Define $\rho _{L,g(w,L)}^{w} (\bar{v} _{w,L})= 1 \setminus \bigcup \{\sigma _{L,i}^{w} (\bar{v} _{w,L}) : i < g(w,L) \}$.  

Part $1$ and $2$ of the conclusion clearly hold.    For $ m < g(w,L),$  $$ \mu (   p_{m}^{w} \triangle \rho _{L,m}^{w} (\bar{v} _{w,L}))\leq \sum _{i \leq m} \mu (  p_{i}^{w} \triangle \sigma _{L,i}^{w} (\bar{v} _{w,L}) $$

$$\leq g(w,L) / 2^{ L+5} L [ g(w,L) ] ^{r^{2}+1}= 1/2^{ L+5} L [ g(w, L) ] ^{r^{2}}. $$

For $m= g(w,L)$  

\begin{eqnarray*}
 & \mu (  p_{g(w,L)}^{w} \triangle \rho _{L,g(w,L)}^{w} (\bar{v} _{w,L})  \\
\leq & \sum _{i \leq g(w,L) } \mu (  p_{i}^{w} \triangle \sigma _{L,i}^{w} (\bar{v} _{w,L})) + \mu ( \bigcup \{ p_{i}^{w} : i > g(w,L) \})  \\
 \leq & g(w,L) /( L2^ { L+5} [g(w,L)]) + 1/L2^{ L+5}.& \\
\end{eqnarray*}
This concludes the proof of lemma \ref{lemm99}. \end{proo} 

\begin{lemm}\labl{lemm4.30} (GCH)  Let $ s < \omega $ and for $ 1\leq i \leq s$ let $ h_{i}: [ \aleph _{2}]^{2} \times \omega \longrightarrow
\omega$.  There exists $ A = \langle \alpha _{i} : i < \omega \rangle \in [\aleph _{2}] ^{\omega} $ and for $ 1 \leq i \leq s$ there exist functions $ \hat{h}_{i} : \omega \longrightarrow \omega $ such that $$ \forall n < \omega \forall m \leq n
\forall w \in [ \{ \alpha _{i} : n <i < \omega \}] ^{2}
(h_{i}(w,m)= \hat{h} _{i} (m)).$$

\end{lemm}

\begin{proo}  A standard ramification argument will show that there exists $ Z_{0} \subseteq \aleph_{2} $ of order type $
\aleph_{1} $ such that for $ \alpha< \beta < \gamma $ in $ Z_{0}, L< \omega,
$ and $ 1 \leq i \leq s
(h_{i}(\{ \alpha, \beta\},L ) = h_{i}( \{\alpha, \gamma\} , L ))$. See \cite{Nash, EHMR} for
details.  For $ \alpha \in Z_{0}, L < \omega $ and $ 1 \leq i \leq s $ define $ h_{i,\alpha } ( L )= h_{i}
( \{\alpha, \beta\}, L ) $ where $\beta > \alpha $ is chosen  in $ Z_{0}$.
By cardinality considerations there exists a sequence $ \langle Z_{i} :1 \leq i < \omega \rangle$
of subsets of $Z_{0} $ such that for all $ k< \omega$, we have $ Z_{k+1} \subseteq Z_{k},
 | Z_{k} | = \aleph_{1}$ and for all $ \alpha, \beta \in Z_{k+1}, h_{i,\alpha} \restriction (k+1) = 
h_{i,\beta} \restriction (k+1)$.  We define $A = \{ \alpha _{i} : i  < \omega \} $ in the 
following manner. Let $ \alpha _{0} $ be minimal in $ Z_{1}$ and inductively
define $ \alpha _{i}$ to be minimal in $ Z_{i+1} \setminus \{ \alpha _{0 } , \ldots , \alpha _{i-1} \} $.  We then define the functions $ \hat{h}_{i} $ by $ \hat{h}_{i}(k) = h_{i,\alpha _{k}} (k).$  

To verify the lemma let $ n < \omega $ and $ m \leq n$.  Choose $ w = \{ \alpha _{t}, \alpha _{v} \}
\in [ \{ \alpha _{k} : n < k < \omega \}]^{2}.$  Then for $1 \leq i \leq s \,(
h_{i}(w,m) = h_{i} ( \{ \alpha _{t}, \alpha_{v} \},m)= h_{i, \alpha _{t}}(m)=
h_{i, \alpha _{m}}(m) =\hat{h}_{i}(m)$.  Thus the lemma is proved.
\end{proo}

Let $ b,g : [\aleph _{2} ]^{2} \times \omega
\longrightarrow \omega$ be the functions chosen above and $f, \rho ^{w} _{L,m},
\bar{v} ^{w} _{L} $ satisfy the conclusion of lemma \ref{lemm99}.  Let
 $ A = \langle \alpha _{i }: i < \omega \rangle \in
[\aleph _{2}] ^{ \omega}, \hat{b}, \hat{g}, \hat{f} ; \omega \longrightarrow
\omega $ be the set an functions obtained when the lemma \ref{lemm4.30} is applied with $ s = 3$ and $ (h_{1}, h_{2}, h_{3} ) = ( b,g,f)$.   We now verify that $ [ I, n,n,\hat{g},\hat{f}]$ holds for all
$ n < \omega $.  To this end fix $ n < \omega$.  Define $ t< \omega$
to be $ n + \max \{g(m) : m \leq n \}+1$.  For $ w = \{ i,j \}
\in [ n]^{2}$ define $w^{*}$ to be $\{ \alpha _{t + i}, \alpha _{t + j} \}$.  Then for $ w \in [n]^{2}, L < n, 1 \leq m \leq \hat{g}(L)$ define $ \tau ^{w} _{L,m}$ to
be $ \rho ^{w^{*}} _{L,m}$
and $ \bar{u} _{w,L} $ to be $ \bar{v} _{ w^{*},L}$.

\medskip

We will now verify that C1-C5 hold for these sequences of boolean terms
and generators. C1-C3 will follow from lemma \ref{lemm4.13}, C4 from
lemma \ref{lemm4.15} and C5 from lemma \ref{lemm4.14}.

\medskip
  
\begin{lemm}\labl{lemm4.13} Let $\hat{g}, \hat{f} : \omega \longrightarrow \omega,A\subset \aleph_{2}$
and $ \tau ^{w} _{L,m}, \bar{u} _{w,L}, \, ( w \in
[n]^{2}, L < n, 1 \leq m \leq \hat{g}(L))$ be as defined above.  Then 
\begin{enumerate}
\item $ \bar{u} _{ w,L } = \{ y_{ w,L,k} : k \leq \hat{f}(L) \}$
\item $  \langle \tau _{L,m}^{w} (\bar{u} _{w,L})
: m \leq \hat{g}(L) \rangle $ is a partition sequence
\item  For $ m < \hat{g}(L)$ we have $$ \mu (   p_{m}^{w^{*}} \triangle \tau _{L,m}^{w} (\bar{u} _{w,L})) < \frac{1} {2^{L+3} L[\hat{g}(L)] ^ {r^{2}}  }$$
\item $ \mu (
 p_{\hat{g}(L)}^{w^{*}} \triangle \tau _{L,\hat{g}(L)}^{w} (\bar{u} _{w,L})) < 
\frac{1}{  L2^{L+3} }$.
 
\end{enumerate}
\end{lemm}
\begin{proo}  For $w \in [n]^{2}, L < n \,( g (w^{*} , L) = \hat{g}(L) $ and  $f ( w^{*}, L) = \hat{f}(L))$.  
\end{proo}

\begin{lemm}\labl{lemm4.15}  Let $ w \in [n] ^{2}$ and $N < L < n$.  For the sequences of boolean terms defined above
 $$( \mu (\bigcup \{\tau _{N,m}^{w} (\bar{u} 
_{w,N}) \cap \tau _{L,m}^{w} (\bar{u}_{w,L } ) : m \leq \hat{g}(N ) \} )  > 1 - 1/( 2^{N}).$$
\end{lemm}
\begin{proo}
\begin{eqnarray*}
& \mu (\bigcup \{\tau _{N,m}^{w} (\bar{u} 
_{w,N}) \cap \tau _{L,m}^{w} (\bar{u}_{w,L } ) : m \leq \hat{g}(N ) \} )  \\
\geq & \mu (\bigcup \{\tau _{N,m}^{w} (\bar{u} 
_{w,N}) \cap \tau _{L,m}^{w} (\bar{u}_{w,L })\cap p_{m}^{w^{*}}  : m \leq \hat{g}(N ) \} ) \\
= & 1-[ \mu ((\bigcup \{\tau _{N,m}^{w} (\bar{u} 
_{w,N}) \cap \tau _{L,m} ^{w}(\bar{u}_{w,L })\cap p_{m} ^{w^{*}} : m \leq \hat{g}(N ) \} )^{c} )]  \\
\geq & 1-(  \sum _{n < \hat{g}(N ) } \mu (   p_{n}^{w^{*}} \triangle \tau _{N,n}^{w} (\bar{u} _{w,N}))+\sum _{n < \hat{g}(L) ) } \mu (   p_{n}^{w^{*}} \triangle \tau _{L,n}^{w} (\bar{u} _{w,L}) \\ 
& + \mu ( p_{\hat{g}(N)}^{w^{*}} \triangle \tau _{L,\hat{g}(N)}^{w} (\bar{u} _{w,N}))+ \mu ( p_{\hat{g}(L)}^{w^{*}} \triangle \tau _{L,\hat{g}(L)}^{w} (\bar{u} _{w,L}))+
\mu( \bigcup \{ p_{m }^{w^{*}} :m > \hat{g}(N) \}) \\
 \geq & 1- (3/2^ {N +2 }) \\ 
\geq & 1-1/2^{N}.\\
\end{eqnarray*}
This concludes the proof of lemma \ref{lemm4.15}.
\end{proo}
\begin{lemm}\labl{lemm4.14}  Let $ L < n $ and $P \in [n]^{r}$.  The measure of 
$$ \bigcup \{  \bigcap  \{\tau ^{z}_{L,c(z)}(\bar{u} _{z,L}) : z \in [P]^{2} \} : c \in 
{\cal C}(P,L) \wedge c {\text{ realizes } } I \} $$ is less than $ 1/L$.

\end{lemm}

\begin{proo} First note that for $ z \in [P]^{2}$ and $ 1 \leq m \leq \hat{g}(L),$ 
$$p^{z^{*}}_{m}/ J \Vdash d(z^{*}) = b(z^{*},m).$$  Now $ z ^{*} \in [ \{ \alpha _{s} : s \geq t \} ] ^{2}$ and $ m < t$ so $ b ( z ^{*} , m)
= \hat {b} (m) $.  Thus, for $ c \in {\cal C} (P,L)$, $$ q = _{\defin}
\bigcap \{ p_{c(z)} ^{z^{*}} : z \in [P]^{2} \}/J \Vdash  ( \forall z \in [P]^{2}
(d(z^{*}) = b(c(z))))$$ if $ q \not= J$.  Thus if $ c $ realizes
$I$ on $P$ and $ q \not= J$ then in some generic extension, $d$ realizes $I$ on $P^{*} = \{ \alpha _{t + i} : i \in P \}$.   Since we
assume that $d$ does not realize $I$ we can conclude that $ q=J$ and $ \mu(\bigcap \{ p_{c(z)} ^{z^{*}} : z \in [P]^{2} \}) =0$.  Secondly note that $ | {\cal C} (P,L)| < g(L) ^{r^{2}}.$  

We first examine those colorings that induce $I$ and involve at least one color other than $g(L)$.  For each such $c$, $$ \mu (\bigcap \{ \tau ^{z} _{L,c(z)}(\bar{u}_{z,L}) : 
z \in [P]^{2}  \}) \leq \min \{\mu(\tau_ { L ,c(z)}^{z} 
( \bar{u} _{z,L } ) \triangle p_{  c(z)}^{z^{*}}: z \in [P]^{2}\}.$$  By lemma
\ref{lemm4.13} this measure is at most $ 1/ ( 2L [ g(L)] ^
 {r{^2}}
) $.  Thus the probability of any of the colorings under consideration inducing
$I$ is less than $ 1/2L$.  In the case that the coloring induces $I$ and uses  only the color
 $ g(L) $ (implying that there is only one such coloring),
 $$ \mu (\bigcap \{ \tau ^{z} _{L,g(L)}(\bar{u}_{z,L}) : 
z \in [P]^{2}  \}) \leq \min \{\mu(\tau_ { L ,g(L)}^{z} 
( \bar{u} _{z,L } ) \triangle p_{  g(L)}^{z^{*}}: z \in [P]^{2}\}.$$
By lemma \ref{lemm4.13} this value is less than
$ 1/2L$.  Thus lemma \ref{lemm4.14} is proved.
\end{proo}

This finishes the proof of lemma \ref{lemm3} and concludes the proof of the
main theorem.

For the work in this paper, $ \omega$-colorings were defined as mappings from
pairs of ordinals into $ \omega $.  Clearly this can be generalized
so that they are mappings from $r$-tuples of ordinals into $ \omega$.
The concept of an $r$-identity can then be defined as can the collection
of $r$-identities realized by an $ \omega$-coloring, and the collection
(denoted $ {\cal I }^{r} (\kappa)$) of $r$-identities realized by all
$ \omega$-colorings, $ f :[\kappa]^{2} \longrightarrow \omega$.  We
believe that the results of this paper can be extended to show that
${\cal I } ^{r} ( \aleph _{r} ) = {\cal I } ( \aleph _{ \omega })$.
We also believe that these results can be demonstrated by adding many
Cohen reals.

\end{document}